\documentclass[12pt,epsf,amstex]{article}
\usepackage [dvips]{graphicx}
\usepackage{amsmath}
\usepackage{amssymb}
\usepackage{epsfig}
\usepackage{verbatim}

\addtocounter{secnumdepth}{1}
\usepackage[utf8]{inputenc}   
\usepackage[T1]{fontenc}
\usepackage{ae,aecompl}

\begin{document}
\newcommand{\ab}{b}
\newcommand{\bE}{\bar{E}}
\newcommand{\calF}{{\cal F}}
\newcommand{\calH}{{\cal H}}
\newcommand{\calJ}{{\cal J}}
\newcommand{\calM}{{\cal M}}
\newcommand{\calN}{{\cal N}}
\newcommand{\hP}{\hat{P}}
\newcommand{\hPi}{\hat{\Pi}}
\newcommand{\sumn}{\sum_{n=1}^N}

\newcommand{\bnu}{\bar{\nu}}
\newcommand{\bnuone}{\bar{\nu}_1}
\newcommand{\bnutwo}{\bar{\nu}_2}

\newcommand{\veca}{\vec{a}}
\newcommand{\vecA}{\vec{A}}
\newcommand{\ai}{{a_i}}
\newcommand{\aone}{{a_1}}
\newcommand{\atwo}{{a_2}}
\newcommand{\Ai}{{A_i}}
\newcommand{\Aone}{{A_1}}
\newcommand{\Atwo}{{A_2}}

\newcommand{\vecalpha}{\vec{\alpha}}
\newcommand{\vecg}{\vec{g}}
\newcommand{\vecp}{\vec{p}}

\newcommand{\tA}{\tilde{A}}
\newcommand{\tB}{\tilde{B}}
\newcommand{\tP}{\tilde{P}}
\newcommand{\tbeta}{\tilde{\beta}}
\newcommand{\tgamma}{\tilde{\gamma}}
\newcommand{\tcalM}{\widetilde{\cal M}}
\newcommand{\betast}{{\beta_*}}
\newcommand{\fstar}{f_*}

\newcommand{\intp}{\int_{-\pi}^{\pi}\frac{\dd p}{2\pi}}
\newcommand{\intpone}{\int_{-\pi}^{\pi}\frac{\dd p_1}{2\pi}}
\newcommand{\intptwo}{\int_{-\pi}^{\pi}\frac{\dd p_2}{2\pi}}
\newcommand{\ointz}{\oint\frac{\dd z}{2\pi{\rm i}}}
\newcommand{\qext}{q_{\rm ext}}

\newcommand{\bO}{{\bf{O}}}
\newcommand{\bR}{{\bf{R}}}
\newcommand{\bS}{{\bf{S}}}
\newcommand{\bT}{\mbox{\bf T}}
\newcommand{\bt}{\mbox{\bf t}}
\newcommand{\half}{\frac{1}{2}}
\newcommand{\thalf}{\tfrac{1}{2}}
\newcommand{\bsA}{\mathbf{A}}
\newcommand{\bsV}{\mathbf{V}}
\newcommand{\bsE}{\mathbf{E}}
\newcommand{\bsT}{\mathbf{T}}
\newcommand{\bsZ}{\hat{\mathbf{Z}}}
\newcommand{\bse}{\mbox{\bf{1}}}

\newcommand{\invup}{\rule{0ex}{2ex}}

\newcommand{\bGamma}{\boldmath$\Gamma$\unboldmath}
\newcommand{\dd}{\mbox{d}}
\newcommand{\ee}{\mbox{e}}
\newcommand{\p}{\partial}

\newcommand{\cdottt}{\!\cdot\!}
\newcommand{\wt}{\widetilde}
\newcommand{\wh}{\widehat}

\newcommand{\la}{\langle}
\newcommand{\ra}{\rangle}

\newcommand{\beq}{\begin{equation}}
\newcommand{\eeq}{\end{equation}}
\newcommand{\bea}{\begin{eqnarray}}
\newcommand{\eea}{\end{eqnarray}}
\def\lsim{\:\raisebox{-0.5ex}{$\stackrel{\textstyle<}{\sim}$}\:}
\def\gsim{\:\raisebox{-0.5ex}{$\stackrel{\textstyle>}{\sim}$}\:}

\numberwithin{equation}{section}

\thispagestyle{empty}
\title{\Large {\bf Mixed-strategy Nash equilibrium}\\[2mm] 
{\bf for a discontinuous symmetric $N$-player game}}

\author{{H.J. Hilhorst and C. Appert-Rolland}\\[5mm]
{\small Laboratoire de Physique Th\'eorique (UMR 8627)}\\[-1mm] 
{\small CNRS, Universit\'e Paris-Sud}\\
{\small Universit\'e Paris-Saclay, 91405 Orsay Cedex, France}\\}

\maketitle
\begin{small}
\begin{abstract}
We consider a game in which each player must
find a compromise between more daring strategies
that carry a high risk for him to be eliminated,
and more cautious ones
that, however, reduce his final score.
For two symmetric players this game was
originally formulated in 1961 by Dresher,
who modeled a duel between two opponents.
The game has also been of interest
in the description of athletic competitions. 
We extend here the two-player game to an arbitrary number $N$ of 
symmetric players.
We show that there is a mixed-strategy Nash equilibrium and find
its exact analytic expression,
which we analyze in particular in the limit of large $N$,
where mean-field behavior occurs.
The original game with $N=2$ arises as a singular limit of the general case. 
\end{abstract}
\end{small}
\vspace{12mm}

\newpage 


\section{Introduction}

The interest of physicists in game theory is increasing. 
The role of this theory
in population dynamics, phase transitions, and other traditional areas of
statistical physics is by now well-documented
(see, {\it e.g.} \cite{hauertszabo2005}); and the new field of
quantum games (see {\it e.g.} 
\cite{benjaminhayden2001,piotrowskisladkowski2003,cheonetal2006})
is blossoming.
Moreover, the methods of statistical physics combined with those of game
theory find applications in other areas of science
such as information theory \cite{wolpert2006},
linguistics \cite{dallastaetal2006}, and the social sciences
\cite{degondetal2014}. 

Game theory and 
statistical physics are both probabilistic and their methods
meet whenever there is a large number $N$ of
participating players. The players 
will typically be particles in the case of statistical physics, 
may represent economic agents in an application to economics,
or be still some
other kind of fundamental entity in yet another domain of application.
\vspace{2mm}

In this short paper we consider $N$ equivalent players $j=1,2,\ldots,N$
participating in a game that obeys the following rules.

1. Each player $j$ chooses an $x_j$ in the real interval
$[0,1]$, called his {\it strategy space.}

2. Player $j$ is randomly eliminated from the game with probability $x_j$.

3. Among the players not eliminated, the one having the 
largest $x_j$ is the winner; any tie is broken randomly.
If it so happens that all players are eliminated, the winner is
chosen at random among them.

4. The winner gains a payoff $N-1$, all others lose $1$. Hence this is a
zero-sum game.

A mathematically important feature of this game is
the discontinuity of the payoff as a function of
$(x_1,\ldots,x_N)$: the winner's profit drops dis\-con\-tinuous\-ly
when his score decreases continuously below the score of the runner-up.
\vspace{2mm}

This model was originally formulated for $N=2$ by Dresher \cite{dresher1961} 
and meant to represent a duel between two opponents.
Another interpretation \cite{appertrollandetal2017}
is that of a running competition in which 
runner $j$ invests an ``energy'' $x_j$.
The more energy he invests, the higher his risk of being eliminated
(say, by exhaustion or any other kind of misfortune beyond his control), 
but also the higher his chances of winning if not eliminated.
Of course such simple models do not do justice to the details of any
specific competition; however, they serve to bring out a few general
principles that play a role.

Let us look at, say, player $j$.
A {\it strategy\,} of player $j$ is a probability distribution 
$f_j(x_j)$ on $[0,1]$
from which at each round of the game he chooses an $x_j$ randomly and
independently. If $f_j(x_j)=\delta(x_j-X_j)$ for some $X_j\in[0,1]$,
we say that player $j$ has the {\it pure\,} strategy $X_j$;
in all other cases we call his strategy {\it mixed.}
A set $\big(f^*_1(x_1),\ldots,f^*_N(x_N)\big)$ of $N$ strategies is
called a {\it Nash equilibrium\,}
if none of the players can improve his gain by 
unilaterally deviating from his strategy while all others stick to theirs.

We will show by explicit construction that
the game presented above possesses a symmetric Nash equilibrium,
{\it i.e.,} one in which 
$f^*_j(x)=\fstar(x)$ for all $j=1,2,\ldots,N$.
For short and with a slight abuse of language
we will also refer to $\fstar(x)$ as the Nash equilibrium of the game.
By symmetry we then know that in such a symmetric equilibrium
the expected gain for any one player is zero.
However, finding $\fstar(x)$ is a nontrivial problem
that cannot be solved by symmetry considerations alone.

For $N=2$ this game has become a textbook example. 
Its solution is \cite{dresher1961}
\beq
\fstar(x) = \left\{
\begin{array}{lc}
\dfrac{1}{4(1-x)^3}\,, \phantom{XX}& 0\leq x\leq \tfrac{2}{3}\,,\\[4mm]
0, & \tfrac{2}{3}<x\leq 1.
\end{array}
\right.
\label{solN2}
\eeq
The fact that there should exist an interval where $\fstar(x)=0$ is 
particularly unintuitive. 

In recent work \cite{appertrollandetal2017}
the $N=2$ game was extended to
two {\it asymmetric\,} players, for which the Nash equilibrium is a pair 
$\big( f^*_1(x_1),f^*_2(x_2) \big)$ of distinct functions. 
In this note we extend the Nash equilibrium solution (\ref{solN2}) to the
case of $N>2$ symmetric players.
Such multi-player games typically lead to nonlinear problems and
we are not aware of the existence of any analytic solutions of the kind 
that we present here.
Simplified but exactly solvable models as this one may illustrate 
general principles or simply be of interest for their own sake;
depending on the context 
they may also serve as paradigms or as testing grounds for numerical methods.

Before entering upon the actual calculation
we will in the following two subsections 
discuss two particular aspects of this problem. 


\subsection{Discontinuous games}

The important subclass of symmetric $N$-player games was considered 
as early as 1951 by Nash \cite{nash1951}. 
Since then, mainstream literature, spread across many areas of science, 
has dealt with proving the
existence of Nash equilibria under a diversity of conditions, developing search 
algorithms, and studying their algorithmic complexity
(see {\it e.g.} \cite{brandtetal2009}).
Usually, such work first focuses on games with a discrete strategy
space (see {\it e.g.} \cite{hefti2017}), 
continuous strategy spaces such as the interval $[0,1]$ being much
harder to study.
The difficulty still increases considerably when, as is the case here,
the payoff function has discontinuities,
in which case one speaks of a ``discontinuous game.''

Under certain broad conditions \cite{chengetal2004}, 
among which continuity of the payoff function,
symmetric $N$-player games with continuous strategy spaces 
are known to admit symmetric {\it pure\,}-strategy Nash equilibria.
This existence is no longer guaranteed when the payoff function has a
discontinuity. In fact, we know that for $N=2$
the game studied in this work has no pure-strategy Nash equilibrium. 
Worse, some symmetric games fail to have symmetric equilibria altogether,
although they may have asymmetric ones, as was shown \cite{fey2012}
for two-player games. 
Several of these existence results appear to carry over \cite{plan2017} to
games with lower than full permutational $N$-player symmetry. 

It is easy to see that the game studied here cannot have
a symmetric pure-strategy Nash equilibrium. Suppose that 
$\fstar(x)=\delta(x-X_*)$ were one. 
If all players played this strategy, 
then each of them would have 
an expected gain equal to zero.
However, if player 1 were to 
shift his strategy from $x_1=X_*$
to $x_1=X_*+\epsilon$ for an arbitrarily small $\epsilon>0$,
then his probability of being eliminated would increase 
negligibly whereas if not eliminated he would be sure to win.
Hence he would increase his expected gain, contrary to what was supposed.


\subsection{Large-$N$ limit}

In statistical physics a mean-field theory for an $N$-component system
is one in which in the large-$N$ limit
each of the $N$ components (atoms, spins, \ldots)
is subject only to a suitably defined
{\it average\,} effect of the other components.
Depending on the model, mean-field equations may sometimes be derived
exactly, and sometimes require approximations to be made. 
Mean-field theory has often been the first, and sometimes the only possible,
method to answer new questions in physics.
From a statistical physics point of view, 
one may expect -- and we will confirm this below --
that for large $N$ a symmetric $N$-player
game can be described by mean-field equations.
\vspace{2mm}

Our work may be seen
against the background of what are commonly called ``mean-field games'' 
(MFG), even though it does not itself belong to that class.
MFG games are based on mean-field ideas from physics 
that were initially brought to bear on game theory
by Lasry and Lions \cite{lasrylions2006} and by
Huang {\it et al.} \cite{huangetal2006}. These authors 
considered a set of $N$ time evolution equations in which  
the strategy of each player is a time-dependent control function.
Each player's payoff functional depends continuously on the control
function and the system has a pure-strategy Nash equilibrium.

Recent additions to this MFG class of games 
are due to Degond {\it et al.}\,\cite{degondetal2014},
who aim at an application to the social sciences,   
and to Swiecicki {\it et al.}\,\cite{swiecickietal2016} 
and Ullmo {\it et al.}\,\cite{ullmoetal2017},
who establish a connection 
between MFG and various equations of physics.
Explicit solutions of certain specific MFG were given 
by Bardi \cite{bardi2012}. 
The convergence, as $N\to\infty$, 
of the Nash equilibrium of a symmetric $N$-player game to 
the solution of the limiting MFG has been
addressed by Fischer \cite{fischer2017}. 
\vspace{2mm}

We may, perhaps, view the results of this note as a first step towards
constructing the solution of a MFG with {\it discontinuous\,} payoffs. 
These include games where the winner from among a large number $N$ of
competitors takes all the profit,
as does the gold medalist in, for example, a multi-runner marathon.


\section{Payoff in a single round and expected gain}


\subsection{Payoff functions  $G_j(x_1,x_2,\ldots,x_N)$}

\noindent 
In a single round of the game, 
let the $N$ players choose strategies $x_1,x_2,\ldots,$$x_N$
from given distributions $f_1(x_1),\ldots,f_N(x_N)$.
Consider an arbitrary player, say player 1, and write
$G_1(x_1,x_2,\ldots,x_N)$ for his payoff 
averaged over the random elimination process. We then have
\bea
G_1(x_1,x_2,\ldots,x_N) &=& P(x_1,x_2,\ldots,x_N)\times (N-1) 
+ Q(x_1,x_2,\ldots,x_N)\times (-1)
\nonumber\\[2mm]
&=& N P(x_1,x_2,\ldots,x_N) - 1,
\label{xG1}
\eea
in which $P$ and $Q=1-P$ are the 
elimination averaged probabilities for player 1 to win and to
lose, leading to payoffs $N-1$ and $-1$, respectively. 

In order to establish the explicit expression for the probability
$P(x_1,x_2,$$\ldots,x_N)$ that player 1 be the winner given 
$x_1,x_2,\ldots,x_N$,
we must sum over all different subgroups of noneliminated players.
For $n=0,1,\ldots,N-1$, let $J_n$ be a subset of $n$ elements of the
set of indices $\{2,3,\ldots,N\}$.
We then have
\bea
P(x_1,x_2,\ldots,x_N) &=& (1-x_1)\sum_{n=0}^{N-1} \sum_{J_n}
\left[\prod_{j\in J_n}(1-x_j)\theta(x_1-x_j)\right]
\left[\prod_{\substack{j \notin J_n \\ j>1}} x_j\right] 
\nonumber\\[2mm]
&& +\, \frac{1}{N}\,x_1\prod_{j=2}^N x_j\,,
\label{xP}
\eea
where we define the Heaviside step function by
$\theta(x)=0$ for $x\leq 0$ and $\theta(x)=1$ for $x>0$. In equation
(\ref{xP}) the first term on the RHS accounts for all different ways 
for player 1 to win without being eliminated,
and the second one represents his probability of winning
when he himself, as well as all the other $N-1$ players, are eliminated.
For $n=0$ the only possibility for $J_0$ is to be the empty set
and the product on $j\in J_0$ is equal to unity; similarly, for $n=N-1$
the product on $j\notin J_{N-1}$ equals unity.

As is clear from the rules of the game stated in the introduction,
equation (\ref{xP}) for $P$ is valid only when $x_1$ 
is nondegenerate with any of the other $x_j$.
The expression for $P$ in the degenerate case, although easily written down, 
will not be needed if we decide to limit our considerations to
mixed strategies $f_j(x_j)$ that
are sufficiently smooth functions of their argument; in that case 
degeneracies occur with probability zero and do not contribute
to any of the calculations. 
\vspace{2mm}

When (\ref{xP}) is substituted in (\ref{xG1}) we obtain the payoff function
for player 1. The payoff $G_j(x_1,x_2,\ldots,x_N)$ for the $j$th player 
results from a simple permutation of indices. 
The functions $G_j$ fully define the game.

If we specialize to $N=2$ we obtain from (\ref{xP}) and (\ref{xG1}) that  
\beq
G_1(x_1,x_2) = \left\{
\begin{array}{ll}
-1 + 2x_2 - x_1x_2\,, \phantom{XX} & x_1<x_2\,, \\[2mm]
\phantom{-}1 - 2x_1 + x_1x_2\,,     & x_1>x_2\,, \qquad N=2,
\end{array}
\right. 
\label{xGN2}
\eeq
which is the example due to Dresher \cite{dresher1961}.

 
\subsection{Expected gain $\overline{G}_1[x_1;f]$}

Suppose now that the players $2,3,\ldots,N$ all adopt the same 
(sufficiently smooth, but otherwise arbitrary) strategy
$f(x)$ and that player 1 chooses a specific $x_1\in[0,1]$.
Let $\overline{G}_1[x_1;f]$ denote player 1's expected gain when the game is
repeated over many rounds under these
circumstances, that is, the overbar denotes the average with respect to the
strategies $f(x_2), f(x_3),\ldots,f(x_N)$ of the other $N-1$ players.
Using (\ref{xG1}) we may also write
\beq
\overline{G}_1[x_1;f] =
N \overline{P}[x_1;f] - 1,
\label{xPav}
\eeq
in which
\bea
\overline{P}[x_1;f] &=& \int\dd x_2\,f(x_2)...\int\dd x_N\,f(x_N)
P(x_1,x_2,\ldots,x_N) \nonumber\\[2mm]
&=& (1-x_1)\sum_{n=0}^{N-1}\sum_{J_n}
\left[\int_0^{x_1}\dd x\,(1-x)f(x) \right]^n
\left[\int_0^1\dd x\,xf(x)\right]^{N-1-n} \nonumber\\[2mm]
&&\,+\,\frac{1}{N}\,x_1\left[\int_0^1\dd x\,xf(x)\right]^{N-1} 
\label{xPav1}\nonumber
\eea
\bea
\phantom{\overline{P}[x_1;f]} &=& (1-x_1)
\left[\int_0^{x_1}\dd x\,(1-x)f(x) +\int_0^1\dd x\,xf(x)\right]^{N-1}
\,+\,\frac{1}{N}\,x_1\left[\int_0^1\dd x\,xf(x)\right]^{N-1}. 
\nonumber\\[2mm] 
\label{xPav2}
\eea
Here, in going from the first to the second line, we
employed (\ref{xP}), and to go from the second to the third one we used that
the sum on $J_n$ simply leads to an extra factor $\binom{N-1}{n}$.  
Equation (\ref{xPav2}) brings out that
the strategy $x_1$ of player 1 is coupled to two $(N-1)$th powers
representing the $N-1$ other players, but in which their individuality has
disappeared. This property will ensure straightforwardly
that for $N\to\infty$ this game has a mean-field limit.  

We rewrite (\ref{xPav2}) as
\beq
\overline{P}[x;f] = (1-x)U^{N-1}(x) + \frac{1}{N}\,x\,V^{N-1}\,,
\qquad 0\leq x\leq 1,
\label{xPUV}
\eeq
where we introduced the abbreviations
\bea
U(x) &=& \int_0^{x}\dd s\,(1-s)f(s) +\int_0^1\dd s\,sf(s), 
\quad 0\leq x\leq 1,
\label{xU} \\[2mm]
V \equiv U(0) &=& \int_0^1\dd s\,sf(s).
\label{xV}
\eea
The constant $V$ is the average fraction of times 
that a player adopting the mixed strategy 
$f(x)$ is eliminated.
For later reference we observe that
\beq
U^\prime(x) = (1-x)f(x), \qquad 0\leq x\leq 1,
\label{xDU}
\eeq
where the prime denotes differentiation with respect to $x$.


\section{Equation for the Nash equilibrium $\fstar(x)$}


\subsection{Equation for $\fstar(x)$ and its solution}

Game theory tells us that if $\fstar(x)$ is a Nash equilibrium, then
the expected gain $\overline{G}_1[x_1;\fstar]$ is equal to some constant $K$ 
for all $x_1$ in the support of $\fstar$ (written `supp $\fstar$'), 
{\it i.e.} all $x_1$ having $\fstar(x_1)>0$.
If such were not the case, player 1 could increase his expected gain by
putting more weight on the values of $x_1$ leading to a higher payoff.
The symmetry between the $N$ players dictates that 
$\overline{G}_1[x_1;\fstar] = K = 0$,
and consequently, by (\ref{xPav}), we have that 
$\overline{P}[x_1;\fstar]=N^{-1}$ for $x_1\in$\,supp\,$\fstar$. 

We must be prepared for the possibility that supp\,$\fstar$
is only a subinterval of $[0,1]$.
Let us assume, therefore, that there exists an interval on which $\fstar(x)$
is positive, continuous, and sufficiently differentiable.
Then, knowing the value $1/N$ of $\overline{P}[x_1;\fstar]$ on this interval, 
we must have 
\beq
(1-x)U_*^{N-1}(x) + \frac{1}{N}\,x\,V_*^{N-1} = \frac{1}{N}\,,
\qquad x\in\mbox{supp}\,\fstar\,,
\label{A0}
\eeq
in which $U_*$ and $V_*$ are given by (\ref{xU}) and
(\ref{xV}) but with  $\fstar$ instead of $f$.
We will attempt to find $\fstar(x)$ by the method \cite{dresher1961} that 
worked for the case $N=2$. Differentiating (\ref{A0}) once with respect to $x$ 
and using (\ref{xDU}) gives
\beq
-U_*^{N-1}(x) + (N-1)(1-x)^2 \fstar(x) U_*^{N-2}(x) 
              + \frac{1}{N}\,V_*^{N-1} = 0,
\label{A1}
\eeq
Differentiating once more and dividing by $(N-1)(1-x)U_*^{N-3}(x)$ yields
\beq
\big[ -3\fstar(x) + (1-x)\fstar^\prime(x) \big]U_*(x) 
      + (N-2)(1-x)^2\fstar^2(x) = 0.
\label{A2}
\eeq
For $N=2$ the second term on the LHS of (\ref{A2}) 
is absent and we have a linear first-order ordinary
differential equation (ODE) for $\fstar(x)$ that is easily 
solved and in the end produces the well-known result (\ref{solN2}).

From here on we will consider the more complicated case $N>2$.
We proceed by solving $U_*(x)$ from (\ref{A2})
in terms of $\fstar(x)$ and $\fstar^\prime(x)$,
differentiate the resulting equality once more, and use (\ref{xDU}).
This leads to
\beq
(1-x)\fstar(x) = (N-2)\frac{\dd}{\dd x}\,\, \frac{(1-x)^2\fstar^2(x)}
{3\fstar(x)-(1-x)\fstar^\prime(x)}\,,
\label{A3a}
\eeq
which amounts to a nonlinear second-order ODE
for $\fstar(x)$. 
Equation (\ref{A3a})
may be cast in several other forms, but we will not do so here.
We just observe that its homogeneity in $\fstar$ facilitates its analysis.
One may check by substitution that this equation has the solution
\beq
\fstar(x) = \frac{C}{(1-x)^{3-\ab}(1+B-x)^{\ab}}\,, \qquad 
\ab\equiv\frac{N-2}{N-1}\,,
\label{xsolf}
\eeq
which depends on the two constants of integration $C>0$ and $B$.
Equations (\ref{A1}) through (\ref{xsolf}), 
which have been derived from (\ref{A0}), are valid for 
$x\in$\,supp\,$\fstar$.

We now observe that the solution (\ref{xsolf}) is not integrable at $x=1$ and, 
guided by the solution for the case $N=2$, take this as an
indication that the support of $\fstar$ is an interval $[0,a]$ with $a<1$. 
The integrations in (\ref{xU}) and (\ref{xV})
then run effectively until the upper limit $x=a$
and the normalization condition is $U_*(a)=1$.


\subsection{Support of $\fstar(x)$}

The solution $\fstar(x)$ depends on the constants $C$, $B$, and $a$, 
which we will now determine by making sure that equations (\ref{A2}),
(\ref{A1}), and (\ref{A0})  are satisfied.
When equation (\ref{xsolf}) for $\fstar(x)$ together with the explicit
expression for $\ab$ is substituted in (\ref{A2}),
we may solve for $U_*(x)$ and find
\beq
U_*(x) = (N-1)\frac{C}{B}\left( \frac{1+B-x}{1-x} \right)^{\frac{1}{N-1}}.
\label{xsolU1}
\eeq 
Since $U_*(x)>0$, we now see that $B>0$.

The normalization condition $U_*(a)=1$ leads together with (\ref{xsolU1})
to an expression for the normalization constant in terms of $a$ and $B$, 
{\it viz.}
\beq
C = \frac{B}{N-1}\left( \frac{1-a}{1+B-a} \right)^{\frac{1}{N-1}},
\label{xC1}
\eeq
which when combined with (\ref{xsolU1}) gives
\beq
U_*(x) = \left( \frac{1-a}{1-x}     \right)^{\frac{1}{N-1}}
       \left( \frac{1+B-x}{1+B-a} \right)^{\frac{1}{N-1}}
\label{xsolU2}
\eeq
and in particular
\beq
V_*=U_*(0)=\left[ \frac{(1+B)(1-a)}{1+B-a} \right]^{\frac{1}{N-1}}.
\label{xV2}
\eeq
The remaining unknowns are $a$ and $B$.

\begin{figure}[htb]
\begin{center}
\scalebox{.45}
{\includegraphics{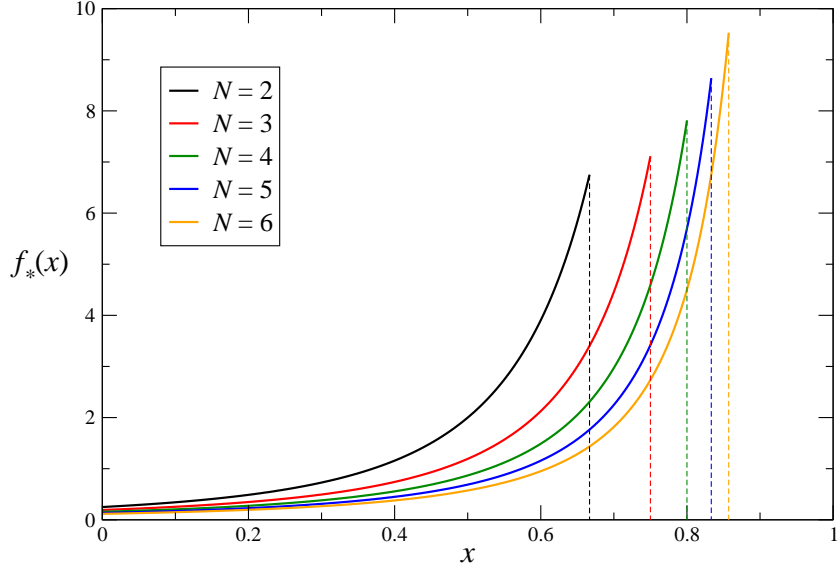}}
\end{center}
\caption{\small 
The Nash equilibrium (\ref{xffinal}) for $N=2,3,4,5,6$.}
\label{figrunN1}
\end{figure}

When (\ref{xsolU2}) and (\ref{xV2}) are substituted in (\ref{A1}) we find
that this equation is satisfied at the condition that
\beq
B = N-1.
\label{xsolB}
\eeq
This leaves the support size $a$ as the only undetermined parameter.
We now return to (\ref{A0}), in which we substitute the results 
(\ref{xsolU2}) for $U_*(x)$ and (\ref{xV2}) for $V_*$, 
as well as (\ref{xC1}) and (\ref{xsolB}) for $C$ and $B$.
It then appears after some algebra that we must have
\beq
a = \frac{N}{N+1}\,.
\label{xsola}
\eeq
We remark that we could not have determined the parameter $a$ this way
if in (\ref{A0}) we had not exploited the fact that the expected gain
$K$ vanishes due to symmetry between all players.  
When (\ref{xsolB}) and (\ref{xsola}) are used in (\ref{xC1}) 
and (\ref{xV2}), we find that
\bea
C &=& N^{-\frac{2}{N-1}}, 
\label{xsolC}\\[2mm]
V_* &=& N^{-\frac{1}{N-1}}.
\label{xsolV}
\eea
After substituting expressions (\ref{xsolB})-(\ref{xsolC}) for 
$B, a$, and $C$ in (\ref{xsolf}) we find that
the Nash equilibrium $\fstar(x)$ is given by
\beq
\fstar(x) = \left\{
\begin{array}{ll}
1/\big[\,N^{2-2\ab}(1-x)^{3-\ab}(N-x)^{\ab}\,\big],\phantom{XX}  & 
0 \leq x \leq a, 
\\[4mm]
0, & a<x\leq 1,
\end{array}
\right.
\label{xffinal}
\eeq
with $\ab = \tfrac{N-2}{N-1}$ and $a=\tfrac{N}{N+1}$.
This is the main result of this note.
Equation(\ref{xffinal}) has been derived for $N>2$. Nevertheless, when
in it we set $N=2$, the equation reduces to the Nash equilibrium
(\ref{solN2}) of the two-player game.
Figure \ref{figrunN1} shows the function $\fstar(x)$ given by (\ref{xffinal}) 
for $N=2,3,4,5,6$.  

\begin{figure}[htb]
\begin{center}
\scalebox{.45}
{\includegraphics{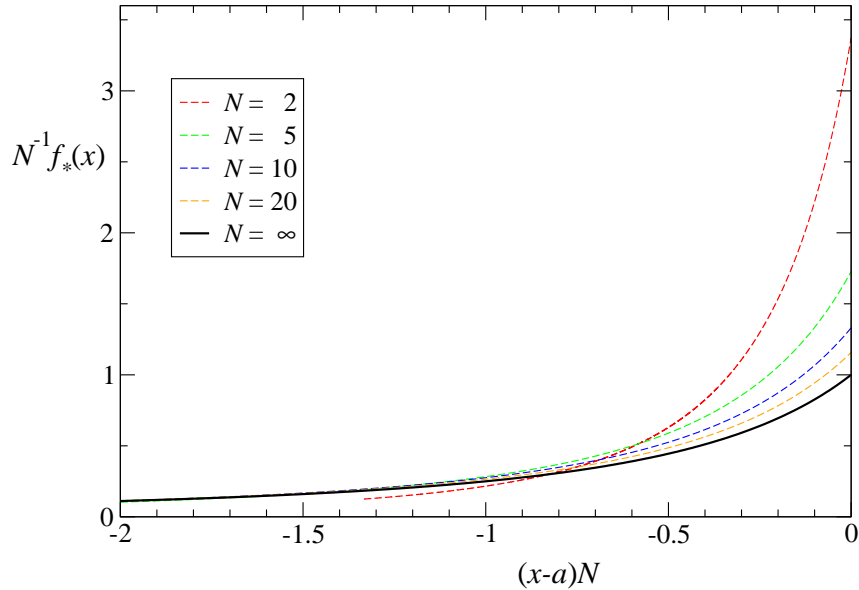}}
\end{center}
\caption{\small The Nash equilibrium shown for
$N=2,5,10,20$ as a function of the scaling variable
$\xi = (x-a)N$. For large $N$ the equilibrium approaches the scaling function
$\calF(\xi)$ (heavy black curve marked $N=\infty$) given by equation
(\ref{xcalF}).}
\label{figrunN3}
\end{figure}


\subsection{Completing the proof}

In order to prove that this is a true Nash equilibrium we have to check one
more thing.
We should verify that when player 1 chooses a strategy 
$x_1$ in the interval $a<x_1\leq 1$, 
while all other players $j=2,3,\ldots,N$ stick to the strategy $\fstar(x_j)$,
then the expected gain of player 1 is less than 0. That is, we have to show
that 
\beq
\overline{P}[x_1;\fstar] < \frac{1}{N}\,, \qquad a<x_1\leq 1.
\label{proofNE}
\eeq
To do so we return to expression (\ref{xPUV}) for $\overline{P}[x_1;\fstar]$.
We observe that for $x_1>a$ the two integrations in (\ref{xU}) 
run until $x=a$, so that $U_*(x_1)=U_*(a)=1$.
Then (\ref{xPUV}) becomes
\bea
\overline{P}[x_1;\fstar] &=& 
(1-x_1)\times 1 + \frac{1}{N}\,x_1\times\frac{1}{N}
\nonumber\\[2mm]
&=&1 - (1-N^{-2})x_1, 
\nonumber\\[2mm]
&<& \frac{1}{N}  \qquad \mbox{ for } \quad \tfrac{N}{N+1}<x_1\leq 1,
\eea
which is the final step in proving that we found a Nash equilibrium.


\subsection{Scaling in the large-$N$ limit}
\label{secscalinglimit}


\subsubsection{Scaling of the distribution}
\label{secscalingdistr}

We study the limit of large $N$ in greater detail.
It appears that as $N\to\infty$ the Nash equilibrium $\fstar(x)$ has almost 
all its weight concentrated in a region of width $\sim N^{-1}$ just below $x=a$.
We therefore consider the scaling limit
\beq
N\to\infty, \qquad x-a\to 0, \qquad \xi\equiv(x-a)N  \ \ \mbox{ fixed},
\label{dsclim}
\eeq
in which $\xi$ is nonpositive.
An elementary calculation shows that in this limit 
\beq
\lim_{\rm scaling} N^{-1}\fstar(x)\dd x \equiv \calF(\xi)\,\dd\xi
\eeq 
with
\beq
\calF(\xi)=\frac{1}{(1-\xi)^2}\,, \qquad -\infty < \xi\leq 0,
\label{xcalF}
\eeq
which satisfies the normalization 
$\int_{-\infty}^0 \dd\xi\, \calF(\xi) = 1$.
Figure \ref{figrunN3} shows how the Nash equilibria tend to their scaling
limit as $N$ gets large. The function $\calF(\xi)=1/(1-\xi)^2$ rightfully
deserves the name ``mean-field strategy.''


\subsubsection{Scaling of the average}
\label{secscalingav}

The  average $\langle x\rangle$ of $x$ represents
the average fraction of times that a player is eliminated.
We are interested in the large-$N$ limit of
this quantity with respect to the Nash equilibrium distribution. 

Although we have the exact relation 
$\langle x\rangle = a + \langle\xi\rangle/N$,
the average of $\xi$ with respect to $\calF(\xi)$ does not exist
and we need to reconsider how to take the large-$N$ limit.
The fast way to find the answer is to use (\ref{xV}) and
to observe that in the Nash equilibrium we have $\langle x\rangle = V_*$
for all $N$, with $V_*$ given by (\ref{xsolV}).
We therefore deduce the large-$N$ expansion
\beq
\langle x\rangle = V_* = 1\,-\,\frac{\log N}{N} 
    + {\cal O}\left( \frac{\log^2 N}{N^2}\right), \qquad N\to\infty.
\label{expV}
\eeq
This brings out quantitatively a lesson that
was qualitatively intuitive: as the number $N$ of players increases,
each player must take higher risks to maximize his chances to win. He does so
by choosing strategies closer and closer to unity.
Since the elimination process acts on each player independently,
in a single round of the 
game there will be typically only 
$\big( 1-\langle x\rangle \big)N \sim\!\log N$ noneliminated players. 

Although true,
it would clearly be an unhelpful oversimplification
to state that for $N\to\infty$ the 
Nash equilibrium converges toward a Dirac delta function at $x=1$.
What matters is the approach to this limit,
and we have shown in this section that that approach is nontrivial.


\section{Conclusion}
\label{secconclusion}

We have in this work generalized one of the simplest discontinuous 
two-player games to an arbitrary number $N$ of symmetric players, 
each having a strategy space $[0,1]$.
We showed by construction that this game has
a symmetric Nash equilibrium $f(x)$ on $[0,1]$ and we presented its 
fully explicit analytic expression for all finite $N=2,3,4,\ldots$. 
We found that in the limit $N\to\infty$
the Nash equilibrium takes a scaling form that depends
only on a single combination $\xi$ of $x$ and $N$.
As had to be expected, when the number $N$ of players increases
the Nash equilibrium shifts to higher values of $x$. Nevertheless,
in the general case, just as for $N=2$, there is
an interval $(1-\frac{1}{N+1},1]$
in which no player can profitably choose his strategy.

The original $N=2$ solution due to Dresher \cite{dresher1961}
appears in this work as a singular limit of the general $N$ case,
\vspace{2mm}

We were motivated by our interest \cite{appertrollandetal2017}
in describing an athletic competition,
{\it e.g.} between runners or cyclists. 
In the present description the athletes do not
interact otherwise than probabilistically.
One of our future goals is to incorporate a physical interaction
between the competitors, 
such as the effect of slipstreaming,
and to see if in that case a similar fully analyzable game can be formulated.

In another direction, it would be interesting to describe a competition 
in which each player adapts his strategy in a time-dependent way. 
For a large number $N$ of players this would naturally lead us to 
the problem of formulating a mean-field game with discontinuous payoff. 


\section*{Acknowledgments}

The authors acknowledge discussions with Amandine Aftalion
on the ``runner'' interpretation of this game,
and with Rida Laraki on the literature and the techniques of game theory.


\appendix


\begin{thebibliography}{10}


\bibitem{hauertszabo2005}
C.~Hauert and G.~Szab\'o,
{\it American Journal of Physics\,} 73 (2005) 405 (2005).

\bibitem{benjaminhayden2001}
S.C.~Benjamin and P.M.~Hayden,
{\it Phys. Rev. A\,} 64 (2001) 030301.

\bibitem{piotrowskisladkowski2003}
E.W.~Piotrowski and J.~Sładkowski,  
{\it International Journal of Theoretical Physics\,} 42 (2003) 1089.

\bibitem{cheonetal2006}
T.~Cheon and I.~Tsutsui, 
{\it Phys. Lett. A\,} 348 (2006) 147.

\bibitem{wolpert2006}
D.H.~Wolpert, {\it Information Theory -- The Bridge Connecting Bounded
Rational Game Theory and Statistical Physics.} 
In: D.~Braha, A.~Minai, and Y.~Bar-Yam
(eds.) Complex Engineered Systems. Understanding Complex Systems. Springer,
Berlin, Heidelberg (2006).

\bibitem{dallastaetal2006}
L.~Dall'Asta, A.~Baronchelli, A.~Barrat, and V.~Loreto,
{\it Phys. Rev. E\,} 74 (2006) 036105.

\bibitem{degondetal2014}
P.~Degond, J.-G.~Liu, and C. Ringhofer,
{\it J. Nonlinear Sci.} 24 (2014) 93.

\bibitem{dresher1961}
M.~Dresher, {\it Games of Strategy: Theory and Applications}, 
Prentice-Hall, Inc.,
  1961.

\bibitem{appertrollandetal2017}
C.~Appert-Rolland, H.J.~Hilhorst, and A.~Aftalion,
\texttt{http://arxiv.org/} \texttt{abs/1709.06460}, 
{\it submitted to\,} JSTAT.

\bibitem{nash1951}
J.~Nash, {\it Annals of mathematics\,} 54 (1951) 286.

\bibitem{brandtetal2009}
F.~Brandt, F.~Fischer, and M.~Holzer,
{\it Journal of Computer and System Sciences\,} 75 (2009) 163.

\bibitem{hefti2017}
A.~Hefti, {\it Theoretical Economics\,} 12 (2017) 979.  

\bibitem{chengetal2004}
S.-F.~Cheng, D.M.~Reeves, Y.~Vorobeychik, and M.P.~Wellman,
{\it Proceedings of the
6th International Workshop On Game Theoretic And Decision Theoretic Agents
GTDT 2004} (2004) p.\,71. Research Collection School Of Information Systems.

\bibitem{fey2012}
M.~Fey,
{\it Games and Economic Behavior} 75 (2012) 424.

\bibitem{plan2017}
A.~Plan, {\it preprint\,} (2017).

\bibitem{lasrylions2006}
J.-M.~Lasry and P.-L.~ Lions,
{\it C.R.~Acad.~Sci.~Paris,} Ser.\,I 343 (2006).

\bibitem{huangetal2006}
M.~Huang, R.P.~Malham\'e, and P.E.~Caines,
{\it Communications in Information and Systems\,} 6 (2006) 221.

\bibitem{swiecickietal2016}
I.~Swiecicki, T.~Gobron, and D.~Ullmo,
{\it Phys. Rev. Lett.} 116 (2016) 128701.

\bibitem{ullmoetal2017}
D.~Ullmo, I.~Swiecicki, and T.~Gobron,
\texttt{arXiv:1708.07730}.

\bibitem{bardi2012}
M.~Bardi, {\it Networks and Heterogeneous Media\,} 
7 (2012) 243.

\bibitem{fischer2017}
M. Fischer, {\it Ann. Appl. Probab.} 27 (2017) 757.


\end{thebibliography}
\end{document}